\documentclass{amsart}
\usepackage{color}
\def\R{\mathbb R}
\def\N{\mathbb N}

\def\C{\mathbb C}
\def\D{\mathbb D}
\def\eps{\varepsilon}

\newcommand{\cal}[1]{{\mathcal #1}}

\newtheorem{thm}{Theorem}

\newtheorem{lem}{Lemma}

\theoremstyle{definition}
\newtheorem{definition}{Definition}
\theoremstyle{remark}
\newtheorem{rem}{Remark}

\begin{document}
\title[Heat equation with dynamical boundary...]
{Heat equation with dynamical boundary conditions of
reactive--diffusive type}
\author{Juan Luis V\'azquez}
\address[J. L. V\'azquez]
       {Dpto. Matem\'aticas, Univ. Aut\'{o}noma de Madrid\\
       28049 Madrid SPAIN}
\email{juanluis.vazquez@uam.es}
\author{Enzo Vitillaro}
\address[E.~Vitillaro]
       {Dipartimento di Matematica ed Informatica, Universit\`a di Perugia\\
       Via Vanvitelli,1 06123 Perugia ITALY}
\email{enzo@dmi.unipg.it}

\date{\today}
%\subjclass{}


%\vskip 2cm

\begin{abstract}
This paper deals with the heat equation  posed  in a bounded regular
domain $\Omega$ of $\R^N$ ($N\ge 2$) coupled   with a dynamical
boundary condition of reactive-diffusive type. In particular we
study the problem
$$\begin{cases} u_t-\Delta u=0 \qquad
&\text{in
$(0,\infty)\times\Omega$,}\\
u_t=ku_\nu+l\Delta_\Gamma u\qquad &\text{on
$(0,\infty)\times \Gamma$,}\\
u(0,x)=u_0(x) &
 \text{on $ \Gamma$,}
\end{cases}$$
where $u=u(t,x)$, $t\ge 0$, $x\in\Omega$, $\Gamma=\partial\Omega$,
$\Delta=\Delta_x$ denotes the Laplacian operator with respect to the
space variable,  while $\Delta_\Gamma$ denotes the Laplace--Beltrami
operator on $\Gamma$, $\nu$ is the outward normal to $\Omega$, and
$k$ and $l$ are given real constants, $l>0$. Well--posedness is
proved for data $u_0\in H^1(\Omega)$ such that ${u_0}_{|\Gamma}\in
H^1(\Gamma)$. We also study higher regularity of the solution.
\end{abstract}

\maketitle

\section{\bf Introduction and main results}

\label{intro} \noindent We deal with the evolution problem
consisting in the standard heat equation posed  in a bounded domain,
supplied with  a dynamical (or Wentzell) boundary condition. The
precise problem is
\begin{equation}\label{P}
\begin{cases} u_t-\Delta u=0 \qquad &\text{in
$(0,\infty)\times\Omega$,}\\
u_t=ku_\nu+l\Delta_\Gamma u\qquad &\text{on
$(0,\infty)\times \Gamma$,}\\
u(0,x)=u_0(x) &
 \text{on $ \Omega$.}
\end{cases}\end{equation}
Here  $u=u(t,x)$, $t\ge 0$, $x\in\Omega$, where $\Omega$ is a
$C^\infty$ regular bounded domain of $\R^N$ ($N\ge 2$) and
$\Gamma=\partial\Omega$. The first equation states the law of standard diffusion or heat conduction in $\Omega$, and $\Delta=\Delta_x$ denotes the Laplacian operator with respect to the space variable.
In the boundary equation \eqref{P}${}_2$, the value of $u$ is  assumed to be the
trace of the  function $u$ defined for $x\in\Omega$, $\Delta_\Gamma$ denotes the Laplace--Beltrami operator on $\Gamma$,  $\nu$ is the outward normal to $\Omega$,  and $k\in\R$ and $l>0$ are given constants; the term $ku_\nu$ represents the interaction domain-boundary, while $l\Delta_\Gamma u$ stands for a boundary diffusion.

A number of authors have studied parabolic problems with dynamical
boundary conditions like \eqref{P}${}_2$.  Note that we can replace
$u_t$  by $\Delta u$  in this boundary condition which leads to the
form known as generalized Wentzell boundary condition. The problem
has been mostly studied the case when there is no Laplacian term on
the boundary condition, i.\,e., when $l=0$. In particular, when $k\le
0$ Problem \eqref{P} is well-posed. See \cite{AMP}, \cite{Eng},
\cite{escher1}, \cite{escher3}, \cite{escher2}, \cite{escher4},
\cite{fgg}, \cite{grobbelaar} \cite{hintermann} in the case $k<0$
which represents a dissipative interaction; the
non-interactive case $k=0$ is rather trivial.  However, when $k>0$ we are
in the presence of a reactive interaction and Problem \eqref{P} is
ill-posed, as shown in the recent papers \cite{bandle} and
\cite{vazvit}. See also \cite{bandlereichel} and \cite{vazvit2} for
the related case $k=k(x)$.

 The question we address in this paper is the following one: is the
situation improved by adding to the dynamical boundary condition a
Laplace--Beltrami correction term  with $l>0$?  The interest of such
a correction  both for the modeling of parabolic and hyperbolic
problems  has been recently pointed out in \cite{goldsteingiselle}.
In particular \eqref{P} describes (see \cite[p.
465]{goldsteingiselle}) a heat conduction process in $\Omega$ with a
 heat source on the boundary which can depend on the heat flux around
the boundary and on the heat flux across it.  The case of
dissipative interaction, $k<0$, has been studied in \cite{coclite1},
\cite{coclite2}, \cite{FGGRLLB} (see also \cite{carmugnolo1},
\cite{carmugnolo2}
 and \cite{mugnoloarxiv1}). It turns out from the
quoted papers that Problem \eqref{P} is well-posed in the framework
of $L^p(\Omega)\times L^p(\Gamma)$, $1\le p\le \infty$. This is to
be expected since both terms in the right-hand side of the boundary
condition have the ``favorable sign''. The aim of this paper  is to
solve the system in the reactive case $k>0$, that is in the usually
ill-posed case. The estimates of the quoted papers did not allow to
cover this case.

A first step in this study has been performed by the authors of the
present paper in \cite{vazvitLLB},  where we consider the Laplace
equation instead  of the heat equation as domain equation. The
modified problem admits a simple functional framework; the paper
 helped the authors  understand the dynamical boundary
condition \eqref{P}${}_2$ and allowed us to formulate the conjecture
that  turns out to be correct, but the arguments used there do not
work for the heat equation. Indeed, a new estimate is needed to deal
with Problem \eqref{P}, which cannot be  obtained in the framework
of $L^2(\Omega)\times L^2(\Gamma)$.

We want to show that Problem \eqref{P} is well-posed in an appropriate setting. We propose to work in the space
\begin{equation}\label{H}
 H=\{(u,v)\in H^1(\Omega)\times H^1(\Gamma): u_{|\Gamma}=v\},
\end{equation}
where $u_{|\Gamma}$ denotes the trace of $u$ on $\Gamma$,  with the
natural topology inherited by $H^1(\Omega)\times H^1(\Gamma)$. Here
and the sequel, we denote for any $s\in \R$, $H^s(\Omega)$ and
$H^s(\Gamma)$ the Sobolev spaces of complex--valued distributions
respectively on $\Omega$ and $\Gamma$ (see \cite{lionsmagenes1} or
\cite{taylor}). For the sake of simplicity we shall identify, when
useful,  $H$ with is isomorphic counterpart $\{u\in H^1(\Omega):
u_{|\Gamma}\in H^1(\Gamma)\}$ through the identification
$(u,u_{|\Gamma})\mapsto u$, so we shall write, without further
mention, $u\in H$ for functions defined on $\Omega$.

Our main result is the following
\begin{thm}\label{theorem1}
For any  $u_0\in H$ Problem \eqref{P} has  a unique  solution
$u=u(u_0)$ such that
\begin{equation}\label{1.1}
\begin{gathered}
u\in C([0,\infty);H^1(\Omega))\cap C^1((0,\infty);H^1(\Omega))\cap
C((0,\infty);H^3(\Omega))\,,\\
u_{|\Gamma}\in C([0,\infty);H^1(\Gamma))\cap
C^1((0,\infty);H^1(\Gamma))\cap C((0,\infty);H^3(\Gamma))\,.
\end{gathered}
\end{equation}
Moreover,
\begin{multline}\label{1.2} \|\nabla
u(t)\|_{L^2(\Omega)}^2+\|d_\Gamma u_{|\Gamma}(t)\|_{L^2(\Gamma)}^2+
\|u_{|\Gamma}(t)\|_{L^2(\Gamma)}^2\\\le e^{2\lambda_0 t}\|\nabla
u_0\|_{L^2(\Omega)}^2+\|d_\Gamma {u_0}_{|\Gamma}\|_{L^2(\Gamma)}^2+
\|{u_0}_{|\Gamma}\|_{L^2(\Gamma)}^2
\end{multline}
for all $t\ge 0$, where $\lambda_0\ge 0$ is a constant depending on
$\Omega$.
Finally, the family of maps $\{u_0\mapsto u(u_0)(t), t\ge 0\}$
extends to an analytic quasi--contractive semigroup in $H$, and
consequently
\begin{equation}\label{1.3}
u\in C^\infty ((0,\infty)\times\overline{\Omega}).
\end{equation}
\end{thm}

 The solutions are in principle complex-valued but it is clear that for
real-valued data the solution is likewise real-valued. As usual,
more regular solutions are obtained for more regular initial data
satisfying usual compatibility conditions. This is the content of
the following regularity result.

\begin{thm}\label{theorem2} If $u_0\in H^{2n+1}(\Omega)$ and ${u_0}_{|\Gamma}\in
H^{2n+1}(\Gamma)$ for some $n\in\N$, and
\begin{equation}\label{1.4}
{(\Delta^i u_0)}_{|\Gamma}=k (\Delta ^{i-1}u_0)_{\nu}+l\Delta_\Gamma
(({\Delta ^{i-1}u_0)}_{|\Gamma}),\qquad\text{for all
$i=1,\ldots,n$,}
\end{equation}
then
\begin{equation}\label{1.5}
%\left\{
\begin{gathered}
u\in C([0,\infty);H^{2n+1}(\Omega))\cap
C^1([0,\infty);H^{2n-1}(\Omega))\cap\ldots\cap C^n([0,\infty); H^1(\Omega)),\\
u_{|\Gamma}\in C([0,\infty);H^{2n+1}(\Gamma))\cap
C^1([0,\infty);H^{2n-1}(\Gamma))\cap\ldots\cap C^n([0,\infty);
H^1(\Gamma)).
\end{gathered}%\right.
\end{equation}
Finally, if $u_0\in C^\infty(\overline{\Omega})$ and \eqref{1.4}
hold for all $i\in\N$, then
\begin{equation}\label{1.6}
u\in C^\infty ([0,\infty)\times\overline{\Omega}).
\end{equation}
\end{thm}

The proofs of Theorems~\ref{theorem1} and \ref{theorem2} rely on the
study of the resolvent problem with  eigenvalue--dependent boundary
condition, that is
\begin{equation}\label{Plambda}
\begin{cases} -\Delta u+\lambda u=h \qquad &\text{in
$\Omega$,}\\
-ku_\nu-l\Delta_\Gamma u+\lambda u=h\qquad &\text{on $\Gamma$,.}
\end{cases}\end{equation}
where $\lambda\in \C$ and $h\in H$. Such type of problems has been studied by some
authors, starting from the classical papers (see \cite{ES1},
\cite{ES2}) to more recent  ones (see \cite{BHLN} and the
bibliography therein). Our result concerning Problem \eqref{Plambda}
is Theorem~\ref{theorem3} below. Finally, we study the limit
behavior of the solution $u$ when $l\to 0^+$ (vanishing boundary
dissipation). See Theorem~\ref{theorem6} below.

The paper is organized as follows. In Section~\ref{section 2}  we
recall some well--known facts and we state some preliminaries. In
Section~\ref{section 3} we analyze the elliptic problem
\eqref{Plambda}, while in Section~\ref{section 4} we apply the
results obtained to Problem \eqref{P}. In Section~\ref{section 5} we
analyze the limit behavior when $l\to 0^+$, while  the final
section contains some comments on future developments.

%%%%%%%%%%%%%%%%%%%%%%%%%%%%%%%%%%%%%%%%%%%%%%%%%%%%%%%%%%%
\section{\bf Preliminaries and functional setting}
\label{section 2} \noindent {\bf Notation.} We  denote by
$\|\cdot\|_p$, $1\le p\le \infty$,  the norm in $L^p(\Omega)$ and,
also the norm in $L^p(\Omega;\R^N)$ since no confusion is expected.  We
denote by $\|\cdot\|_{p,\Gamma}$ the norm in $L^p(\Gamma)$ and
also, when $p=2$,  the $L^2$ norm for square integrable $1$--forms
on $\Gamma$.

\noindent {\bf Laplace--Beltrami operator.} We recall here, for the
reader's convenience, some well--known facts on the
Laplace--Beltrami operator $\Delta_\Gamma$. We refer to \cite{jost}
or \cite{taylor} for more details and proofs. We start by fixing
some notation. Clearly, $\Gamma$ is a Riemannian manifold endowed
with the natural metric inherited from $\R^N$, given in local
coordinates by $(g_{ij})_{i,j=1,\ldots,N-1}$. We denote by $dV$ the
natural volume element on $\Gamma$, given in local coordinates by
$\sqrt g \,\,dy_1\ldots dy_{N-1}$, where $g=\operatorname{det}
(g_{ij})$. We denote by $\nabla_\Gamma$ the Riemannian gradient and
by $d_\Gamma$ the total differential on $\Gamma$. We use the
notation $(\cdot,\cdot)$ for the Riemannian inner product of vectors
while $(\cdot|\cdot)$ is used for the natural scalar product on
$1$-forms on $\Gamma$ associated to the metric. Then, it is clear
that $(d_\Gamma u|d_\Gamma v)=(\nabla_\Gamma u, \nabla_\Gamma v)$
for $u,v\in C^1(\Gamma)$, so the use of vectors or forms in the
sequel is optional.

The Laplace-Beltrami operator $\Delta_\Gamma$ can be at first defined on $C^\infty(\Gamma)$
by the formula
\begin{equation}\label{1}
-\int_\Gamma (\Delta_\Gamma u) \,\bar{v}\,dV=\int _\Gamma (d_\Gamma
u|d_\Gamma v)\,dV
\end{equation}
for any $u,v\in C^\infty(\Gamma)$, and it is given in local
coordinates by
\begin{equation}\label{2}
\Delta_\Gamma u=g^{-1/2}\sum_{i,j=1}^{N-1}\frac\partial{\partial
y_i}\left(g^{ij}g^{1/2}\frac{\partial u}{\partial y_j}\right),
\end{equation}
 where
$(g^{ij})=(g_{ij})^{-1}$ as usual. Clearly, by  \eqref{2}, $\Delta_\Gamma$ can be considered as a
bounded operator from $H^{s+2}(\Gamma)$ to $H^s(\Gamma)$, for any
$s\in\R$. Consequently, Formula \eqref{1} extends by density to
$u,v\in H^1(\Gamma)$, where the integral in the left--hand side has
to be interpreted in the distributional sense, as $\Delta_\Gamma u\in
H^{-1}(\Gamma)$.

\noindent {\bf Remark.} In the sequel, the notation $dV$ will be
dropped from the  boundary integrals; we hope that the reader
will be able to put in the appropriate integration elements in all
formulas.

Since $\Delta_\Gamma 1=0$ the operator is not injective, but  by
\eqref{1} we have
\begin{equation}\label{3}
\int_\Gamma(-\Delta_\Gamma u +u)\bar{u}=\|d_\Gamma
u\|_{L^2(\Gamma)}^2+\|u\|_{L^2(\Gamma)}^2
\end{equation}
so that the operator $L:=-\Delta_\Gamma+1$  is a topological and
algebraic isomorphism between $H^1(\Gamma)$ and $H^{-1}(\Gamma)$.
Moreover, by elliptic regularity (see \cite[p. 309]{taylor}), \
$L^{-1}: H^{k-1}(\Gamma)\to H^{k+1}(\Gamma)$, $k=0,1,2,\ldots$, is
bounded, so \  $L: H^{k+1}(\Gamma)\to H^{k-1}(\Gamma)$ \ is an
isomorphism. By interpolation, $L^{-1}: H^s(\Gamma)\to
H^{s+2}(\Gamma)$ for all $s\in\R$, $s\ge -1$, giving the inverse of
\ $L: H^{s+2}(\Gamma)\to H^s(\Gamma)$. By duality, this fact holds
for all real $s$.

\medskip

 \noindent {\bf Dirichlet--to--Neumann operator}. We  will also need some well-known facts about this operator
that will be used at some technical points.  We
 refer to \cite{lionsmagenes1}  for details and proofs. For any $u\in H^s(\Gamma)$, $s\in\R$, the
non--homogeneous Dirichlet problem
\begin{equation}\label{8}
  \begin{cases}
        \Delta v=0,& \text{in $\Omega$}, \\
        v =u & \text{on $\Gamma$,}
  \end{cases}
\end{equation}
has a unique solution $v\in H^{s+1/2}(\Omega)$, here denoted by
$v=\mathbb{D}u$. Moreover $\D$ is a bounded operator from
$H^s(\Gamma)$ to $H^{s+1/2}(\Omega)$ for all real $s$, and $v$  has
a normal derivative $v_\nu\in H^{s-1}(\Gamma)$. The operator
$u\mapsto v_\nu$, known as the Dirichlet--to--Neumann operator, is
bounded from $H^s(\Gamma)$ to $H^{s-1}(\Gamma)$, and it will be
denoted in the sequel by $\cal{A}$. For all $u,v\in
C^\infty(\Gamma)$, integrating by parts twice we have
\begin{equation}\label{8bis}
\int_\Gamma \cal{A}u\, \bar{v}=\int_\Omega\nabla (\D u)\nabla (\D
v)=\int_\Gamma u\cal{A}\bar{v}\end{equation} which, by density,
holds for all $u,v\in H^1(\Gamma)$.

\medskip

 \noindent {\bf Functional setting}. In the sequel we  equip $H^1(\Gamma)$ with the equivalent norm
in \eqref{3}, so we  denote
\begin{equation}\label{4}
(u,v)_{H^1(\Gamma)}=\int_\Gamma u\overline{v}+\int_\Gamma(d_\Gamma
u|d_\Gamma v), \qquad \|u\|_{H^1(\Gamma)}^2=(u,u)_{H^1(\Gamma)}
\end{equation}
for all $u,v\in H^1(\Gamma)$. Moreover, since $-\Delta_\Gamma
+1:H^2(\Gamma)\to L^2(\Gamma)$ is an isomorphism we can equip
$H^2(\Gamma)$ with the equivalent norm
\begin{equation}\label{5}
(u,v)_{H^2(\Gamma)}=\int_\Gamma
u\overline{v}+\int_\Gamma\Delta_\Gamma u\Delta_\Gamma \overline{v},
\qquad \|u\|_{H^2(\Gamma)}^2=(u,u)_{H^2(\Gamma)}
\end{equation}
for all $u,v\in H^2(\Gamma)$. Moreover, we denote as usual
\begin{equation}\label{6}
\|u\|_{H^1(\Omega)}^2= \|u\|_2^2+\|\nabla u\|_2^2
\end{equation}

\noindent {\sc The space $H$.} We now introduce, as anticipated in
the introduction, the space $H$ given in \eqref{H}, which by the
Trace Theorem is a closed subset of $H^1(\Omega)\times H^1(\Gamma)$,
hence a Hilbert space with respect to the scalar product inherited
from $H^1(\Omega)\times H^1(\Gamma)$. For the sake of simplicity, we
shall drop the notation $u_{|\Gamma}$, when clear, so we shall write
$\|u\|_{2,\Gamma}$, $\int_\Gamma u$, and so on,  for elements of
$H$, through the already mentioned identification
$(u,u_{|\Gamma})\mapsto u$. We equip $H$ with an equivalent norm
which simplifies our calculations. This is the content of the
following
\begin{lem}\label{lemma0}
We set, for any $u,v\in H$,
\begin{equation}\label{10}
(u,v)_H=\int_\Omega \nabla u \nabla v +\int_\Gamma (d_\Gamma
u|d_\Gamma v)+\int_\Gamma u\overline{v},\qquad \|u\|_H^2=(u,u)_H.
\end{equation}
Then $\|\cdot\|_H$ is equivalent in $H$ to the standard norm
inherited by $H^1(\Omega)\times H^1(\Gamma)$.
\end{lem}
\begin{proof}
We just have to show that if we drop $\|\cdot\|_2$ in the standard
norm of $H^1(\Omega)\times H^1(\Gamma)$ we get an equivalent norm.
This follows by  a Poincar\'e-type inequality which says (see
\cite[Theorem~4.4.6]{ziemer} in the real valued case, the extension
to the complex--valued one being trivial) that
$$\left\|u-\int_\Gamma u\right\|_{2^*}\le C_1 \|\nabla u\|_2\qquad\text{for all
$u\in H^1(\Omega)$,}$$ where $C_1=C_1(N,\Omega)>0$, $2^*$ is the
Sobolev critical exponent, i.e. $2^*=2N/(N-2)$ when $N\ge 3$, $1\le
2^*<\infty$ when $N=2$. Consequently, since $\Omega$ is bounded
 and $\Gamma$ is compact, we get
\begin{equation}\label{10.1}
\begin{aligned}
\|u\|_2\le &\left\|u-\int_\Gamma u\right\|_2+\left\|\int_\Gamma
u\right\|_2 \le  C_1\|\nabla u\|_2+\lambda_{N}(\Omega) \int_\Gamma
|u|\\ \le & C_2 \left(\|\nabla u\|_2+\|u\|_{2,\Gamma}\right)
\end{aligned}
\end{equation}
where $\lambda_N$ denotes the usual Lebesgue measure in $\R^N$ and
$C_2=C_2(N, \Omega)>0$. This estimate completes the proof.
\end{proof}

\medskip

\noindent {\sc The space $V$.} We need a further space
\begin{equation}\label{11}
V=\{(u,v)\in H^2(\Omega)\times H^2(\Gamma): u_{|\Gamma}=v\}
\end{equation}
which is naturally embedded in $H$, and it is a Hilbert space with
respect to the scalar product and norm inherited from
$H^2(\Omega)\times H^2(\Gamma)$. As before we  equip it with a
suitable scalar product which induces a norm equivalent to that one.

\begin{lem}\label{lemma0bis} If we set, for any $u,v\in V$,
\begin{equation}\label{13}
(u,v)_V=\int_\Omega \Delta u\Delta\overline{v} +\int_\Gamma
\Delta_\Gamma u\Delta _\Gamma \overline{v}+\int_\Gamma
u\overline{v},\qquad \|u\|_V^2=(u,u)_V\,,
\end{equation}
then $\|\cdot\|_V$ is equivalent  in $V$ to the standard norm
inherited by $H^2(\Omega)\times H^2(\Gamma)$.
\end{lem}
\begin{proof}
It simply follows by elliptic regularity estimates. Indeed, for any
$u\in H^2(\Omega)$ we have (see \cite[p. 202]{lionsmagenes1})
$$\|u\|_{H^2(\Omega)}\le C_3 \left(\|\Delta u
\|_2+\|u_{|\Gamma}\|_{H^{3/2}(\Gamma)}\right)$$ and consequently,
since $H^2(\Gamma)$ is continuously embedded in $H^{3/2}(\Gamma)$,
for any $u\in V$ we get
\begin{equation}\label{10.2}
\|u\|_{H^2(\Omega)}\le \|\Delta u\|_2+\|u\|_{H^2(\Gamma)}
\end{equation}
 which by
\eqref{5} completes the proof.
\end{proof}

%%%%%%%%%%%%%%%%%%%%%%%%%%%%%%%%%%%%%%%%%%%%%%%%%%%%%%%%%%%
\section{\bf Elliptic theory}
\label{section 3} \noindent This section is devoted to study the
solvability of the coupled elliptic system \eqref{Plambda} when
$l>0$, $k\in\R$, $\lambda\in\C$ and $h\in H$.

\noindent {\bf Definition.} By a solution of Problem \eqref{Plambda}
we mean a function  $u\in V$ such that \eqref{Plambda}$_1$ holds
true  in $L^2(\Omega)$, while \eqref{Plambda}$_2$ holds true in
$L^2(\Gamma)$.

Space $V$ was just introduced in \eqref{11}.  Before stating the
main result of this section we introduce, for any $s\ge 1$, the
further space
\begin{equation}\label{Hs}
H^s=\{(u,v)\in H^s(\Omega)\times H^s(\Gamma): u_{|\Gamma}=v\}.
\end{equation}
Clearly, being closed in the product space $H^s(\Omega)\times
H^s(\Gamma)$,  $H^s$ is a Hilbert space equipped with the norm
inherited norm, which we denote by $\|\cdot\|_{H^s}$. Moreover, it
is naturally embedded in $H$ and $H^1=H$, $H^2=V$ (more precisely,
$\|\cdot\|_{H^1}$ and $\|\cdot\|_H$ are merely equivalent, like
$\|\cdot\|_{H^2}$ and $\|\cdot\|_V$).

Our result concerning \eqref{Plambda} is the following

\begin{thm}\label{theorem3}
There is a positive constant $\lambda_0$, depending on $l,k,\Omega,N$, such that for
 $\lambda\in \C$, $\text{Re} \lambda
\ge \lambda_0$ and any $h\in H$  Problem \eqref{Plambda} has a
unique solution $u\in V$, which also belongs to $H^3$. Moreover, if
$h\in H^s$ for some $s\ge 1$, then $u\in H^{s+2}$.

\noindent Finally, there is $C_4=C_4(l,k,\Omega,s,\lambda)>0$ such
that
\begin{equation}\label{31xx}
\|u\|_{H^{s+2}}\le C_4 \|u\|_{H^s}\qquad\text{for all $h\in H^s$.}
\end{equation}
\end{thm}

In order to solve elliptic problems via the variational method it
is useful to introduce a sesquilinear form, which leads to weak
solutions. The most natural way to perform this procedure for
Problem \eqref{Plambda} would be to multiply (at least formally) the
equation $-\Delta u+\lambda u=h$ by  a test function
$\overline{\phi}$ and integrate over $\Omega$ to get
$$-\int_\Omega \Delta u\overline{\phi}+\lambda \int_\Omega
u\overline{\phi}=\int_\Omega h\overline{\phi}.$$  Integrating by
parts, when $u$ is regular enough,
$$\int_\Omega \nabla u\nabla v -\int_\Gamma u_\nu
\overline{\phi}+\lambda\int_\Omega u\overline{\phi}=\int_\Omega
h\overline{\phi}.$$ Then, using the boundary equation in
\eqref{Plambda} we get (when $k\not=0$)
$$\int_\Omega \nabla u\nabla v+\frac 1k \int_\Gamma
h\overline{\phi}-\frac\lambda k \int_\Gamma u\overline{\phi}+\frac
lk \int_\Gamma \Delta_{\Gamma} u \overline{\phi}+\lambda \int_\Omega
u\overline{\phi}=\int_\Omega h\overline{\phi}.$$ Finally, by
\eqref{1} we arrive to
\begin{equation}\label{15bis}
\int_\Omega \nabla u\nabla \phi-\frac lk \int_\Gamma (d_\Gamma
u|d_\Gamma \phi)-\frac{\lambda}k \int_\Gamma
u\overline{\phi}+\lambda \int_\Omega u\overline{\phi}=-\frac 1k
\int_\Gamma h\overline{\phi}+\int_\Omega h\overline{\phi}.
\end{equation}

Now, it is easy to check that the sesquilinear form in the
left-hand side of \eqref{15bis} is indefinite in the case $k>0$, so
this procedure does not produce useful estimates. Thus, one has to
look for a positive definite sesquilinear form, at least for
$\text{Re}\lambda$ large enough. This is exactly the content of the
following two lemmas. The first one introduces the sesquilinear form
which turns out to be appropriate.

\begin{lem}\label{lemma1}
Let $h\in H$. Then $u\in V$ solves Problem \eqref{Plambda} if and
only if
\begin{equation}\label{15}
a_\lambda(u,v)=(h,v)_H\qquad\text{for all $v\in V$,}
\end{equation}
where the sesquilinear form $a_\lambda$ on $V$ is defined by the
formula
\begin{equation}\label{16}
\begin{aligned}
a_\lambda(u,v)=&\int_\Omega \Delta u\,\Delta \overline{v}+l
\int_{\Gamma}\Delta_\Gamma u \,\Delta_\Gamma \overline{v}+\lambda
\int_\Omega \nabla u\nabla v +(\lambda+l)\int_\Gamma (d_\Gamma
u|d_\Gamma v)\\
-&l \int_\Gamma \Delta_\Gamma u \,\overline{v_\nu}+k\int_\Gamma
u_\nu \,\Delta_\Gamma \overline{v}-k \int_\Gamma u_\nu\,
\overline{v_\nu}-k \int_\Gamma u_\nu \,\overline{v} +\lambda
\int_\Gamma u\,\overline{v}.
\end{aligned}
\end{equation}
Moreover in this case  $u\in H^3(\Omega)$ and $u_{|\Gamma}\in
H^3(\Gamma)$.
\end{lem}

\begin{proof} It is divided into several steps.

(i) {\bf Claim. } {\sl If $u\in V$ is a solution of
\eqref{Plambda}, then $u\in H^3(\Omega)$ and $u_{|\Gamma}\in
H^3(\Gamma)$. \normalcolor}
 To recognize that our claim is true we use elliptic regularity both on
$\Omega$ and $\Gamma$ as follows. Since $u\in H^2(\Omega)$ we have
$u_\nu\in H^{1/2}(\Gamma)$ by the Trace Theorem. So, being
$h_{|\Gamma}\in H^1(\Gamma)$ and $u_{|\Gamma}\in H^{1/2}(\Gamma)$,
from \eqref{Plambda}$_2$ it follows that $-\Delta_\Gamma u+u\in
H^{1/2}(\Gamma)$, so that using the the isomorphism property of
$-\Delta_\Gamma+1$, we conclude  that $u_{|\Gamma}\in
H^{5/2}(\Gamma)$. Consequently, using elliptic regularity for
nonhomogeneous Dirichlet problems (\cite[p. 203]{lionsmagenes1}) we
obtain by \eqref{Plambda}$_1$ that $u\in H^3(\Omega)$.

From this, and using the Trace Theorem again, we get $u_\nu\in
H^{3/2}(\Omega)$. Using \eqref{Plambda}$_2$ again we then get
$-\Delta_\Gamma u+u\in H^1(\Gamma)$, so as before $u_{|\Gamma}\in
H^{3}(\Gamma)$, completing the proof of our first claim.

\smallskip

\noindent (ii) {\bf Claim. \sl  If $u\in V$ is a
solution of \eqref{Plambda}, then formula  \eqref{15} holds.} By the
first claim we have $\Delta u\in H^1(\Omega)$. Moreover, by
\eqref{Plambda}$_1$ we get
\begin{equation}\label{star}
(\Delta u)_{|\Gamma} =\lambda u_{|\Gamma}-h_{|\Gamma}\in
H^1(\Gamma).
\end{equation}
Consequently, we get that $\Delta u\in H$, so from
\eqref{Plambda}$_1$  we have
\begin{equation}\label{17}
(-\Delta u, v)_H+\lambda (u,v)_H=(h,v)_H\qquad\text{for all $v\in
H$.}
\end{equation}
Formula \eqref{17} can be written more explicitly, using \eqref{10},
as
\begin{multline}\label{18}
\int_\Omega \nabla (-\Delta u)\nabla v+\int_\Gamma (d_\Gamma
(-\Delta u)|d_\Gamma v)-\int_\Gamma \Delta u \overline{v}\\
+\lambda \int_\Omega \nabla u\nabla v+\lambda\int_\Gamma (d_\Gamma
u|d_\Gamma v)+\lambda \int_\Gamma u\overline{v}=(h,v)_H.
\end{multline}
Now, using \eqref{star} we can write  \eqref{Plambda}$_2$ in the
form
\begin{equation}\label{ball}
(\Delta u)_{|\Gamma}=ku_\nu+l\Delta_\Gamma u_{|\Gamma}.
\end{equation}
Plugging \eqref{ball} into \eqref{18} we get
\begin{multline}\label{19}
\int_\Omega \nabla (-\Delta u)\nabla v-k\int_\Gamma (d_\Gamma
u_\nu|d_\Gamma v)-l\int_\Gamma (d_\Gamma \Delta_\Gamma u)|d_\Gamma
v)-k\int_\Gamma u_\nu \overline{v}\\-l\int_\Gamma \Delta_\Gamma  u
\overline{v} +\lambda \int_\Omega \nabla u\nabla
v+\lambda\int_\Gamma (d_\Gamma u|d_\Gamma v)+\lambda \int_\Gamma
u\overline{v}=(h,v)_H
\end{multline}
for all $v\in H$. Now we restrict to test functions $v\in V$, we
integrate by parts the first integral in \eqref{19} and we use
\eqref{1} in the first one to get
\begin{multline}\label{20}
\int_\Omega \Delta u\Delta \overline{v}-\int_\Gamma \Delta u
\overline{v_\nu}-k\int_\Gamma (d_\Gamma u_\nu|d_\Gamma
v)+l\int_\Gamma \Delta_\Gamma u\Delta_\Gamma\overline{
v}-k\int_\Gamma u_\nu \overline{v}\\-l\int_\Gamma \Delta_\Gamma  u
\overline{v} +\lambda \int_\Omega \nabla u\nabla
v+\lambda\int_\Gamma (d_\Gamma u|d_\Gamma v)+\lambda \int_\Gamma
u\overline{v}=(h,v)_H.
\end{multline}
Plugging  \eqref{ball} once again  in the second integral in the
left--hand side of \eqref{20} and \eqref{1} in the third and sixth
ones we finally get \eqref{15}.

\smallskip

\noindent (iii) To complete the proof, we now suppose that
\eqref{15} holds for some $u\in V$. We have to prove that $u$ solves
\eqref{Plambda}. Integrating by parts in the third integral in
\eqref{16} and in the first one in \eqref{10} we can then write
\eqref{15} as
\begin{equation}\label{22}
\begin{split}
\begin{aligned}
&\int_\Omega \Delta u \Delta \overline{v}+l\int_\Gamma \Delta_\Gamma
u \Delta_\Gamma \overline{v}-\lambda \int_\Omega u\Delta
\overline{v}+\lambda \int_\Gamma u\overline{v_\nu}+(\lambda+l)
\int_\Gamma (d_\Gamma u|d_\Gamma v)\\
-l&\int_\Gamma \Delta_\Gamma u\overline{v_\nu}+k\int_\Gamma u_\nu
\Delta_\Gamma \overline{v}-k\int_\Gamma u_\nu \overline{v_\nu}-k
\int_\Gamma u_\nu \overline{v}+\lambda \int_\Gamma u\overline{v}
\end{aligned}\\
=-\int _\Omega h\Delta \overline{v}+\int_\Gamma h
\overline{v_\nu}+\int_\Gamma (d_\Gamma h|d_\Gamma v)+\int_\Gamma
h\overline{v}
\end{split}
\end{equation}
for all $v\in V$. Using \eqref{1} we can write \eqref{22} as
\begin{equation}\label{22bis}
\begin{split}
\begin{aligned}
&\int_\Omega \Delta u \Delta \overline{v}+l\int_\Gamma \Delta_\Gamma
u \Delta_\Gamma \overline{v}-\lambda \int_\Omega u\Delta
\overline{v}+\lambda \int_\Gamma u\overline{v_\nu}-\lambda
\int_\Gamma u\Delta_\Gamma
\overline{v}\\
-l&\int_\Gamma \Delta_\Gamma u\overline{v}-l\int_\Gamma
\Delta_\Gamma u\overline{v_\nu}+k\int_\Gamma u_\nu \Delta_\Gamma
\overline{v}-k\int_\Gamma u_\nu \overline{v_\nu}-k \int_\Gamma u_\nu
\overline{v}+\lambda \int_\Gamma u\overline{v}
\end{aligned}\\
=-\int _\Omega h\Delta \overline{v}+\int_\Gamma h
\overline{v_\nu}-\int_\Gamma h\Delta_\Gamma \overline{v}+\int_\Gamma
h\overline{v},
\end{split}
\end{equation}
that is, by grouping the terms with respect to the test function,
\begin{equation}\label{23.2}
\int_\Omega \left (\Delta u-\lambda u+h\right)\Delta
\overline{v}+\int_\Gamma \left(-l\Delta_\Gamma u-ku_\nu+\lambda u-h
\right)\overline{\left(-\Delta_\Gamma v+v_\nu+v\right)}=0.
\end{equation}

The form of \eqref{23.2} suggests now how to proceed. Indeed if we
restrict to test functions $v\in C^\infty_c(\Omega)$, at least to
get \eqref{Plambda}$_1$, we get that $\Delta (-\Delta u+\lambda
u-h)=0$ in distributional sense, which is not \eqref{Plambda}$_1$.
Then it is more useful to start by proving \eqref{Plambda}$_2$. With
this aim, we restrict \eqref{23.2} to test functions $\mathbb{D}v$,
where  $v\in H^2(\Gamma)$, to get
\begin{equation}\label{23.21}
\int_\Gamma \left(-l\Delta_\Gamma u-ku_\nu+\lambda u-h
\right)\overline{\left(-\Delta_\Gamma v+\cal{A}v+v\right)}=0,
\end{equation}
where $\cal{A}$ denotes the Dirichlet--to--Neumann operator already
introduced.

We now claim that by \eqref{23.21} it follows that
\begin{equation}\label{23.3}
\int_\Gamma \left(-l\Delta_\Gamma u-ku_\nu+\lambda u-h
\right)\overline{\phi}=0\qquad\text{for all $\phi\in L^2(\Gamma)$},
\end{equation}
from which clearly one has that \eqref{Plambda}$_2$ holds in
$L^2(\Gamma)$. To prove our claim it is enough to recognize that,
given an arbitrary $\phi\in L^2(\Gamma)$, the problem
 \begin{equation}\label{++}
 -\Delta_\Gamma v+\cal{A}v+v=\phi
 \end{equation}
has a solution $w\in H^{2}(\Gamma)$, which turns out to be
 unique. Hence, our claim is nothing but a refinement, in this particular
case, of a previous result of the authors \cite[Lemma~1]{vazvitLLB}
which says that given $\tilde{l}>0$ and $\tilde{k}\in \R$ there is
$\overline{\Lambda}\ge 0$ such that for $\Lambda\ge
\overline{\Lambda}$ the problem
\begin{equation}\label{+++}
 -\tilde{l}\Delta_\Gamma v-\tilde{k}\cal{A}v+\Lambda v=\phi,
 \end{equation}
has unique solution $v\in H^2(\Gamma)$ for any $\phi\in
L^2(\Gamma)$. In particular, our claim is proved if we prove that,
when $\tilde{k}<0$, then we can take $\overline{\Lambda}=1$. To
prove this fact we argue as in the quoted paper, writing \eqref{+++}
in the more explicit form
\begin{equation}\label{++++}
-\int_\Gamma \tilde{k}\cal{A}v \overline{\psi}+\tilde{l}\int_\Gamma
(d_\Gamma v|d_\Gamma \psi)+\Lambda \int_\Gamma
v\overline{\psi}=\int_\Gamma \phi\overline{\psi}\qquad\text{for all
$\psi\in H^1(\Gamma)$,}
\end{equation}
and then we apply Lax--Milgram theorem (see \cite[p.
376]{dautraylionsvol2}) to the sesquilinear form
$$
a(v,\psi)=-\int_\Gamma
\tilde{k}\cal{A}v\,\bar{\psi}+\tilde{l}\int_\Gamma (d_\Gamma
v|d_\Gamma \psi)+\Lambda \int_\Gamma v\bar{\psi},\qquad v,\psi\in
H^1(\Gamma),
$$
which is trivially Hermitian (by \eqref{8bis}) and continuous. To
recognize that it is also coercive for $\Lambda\ge 1$ we simplify
the argument of \cite{vazvitLLB}. Indeed, since $\tilde{k}<0$ we
have by \eqref{8bis}
$$a(v,v)=-\tilde{k}\|\nabla (\mathbb{D} v)\|_2^2+\tilde{l}\|d_\Gamma v\|_{2,\Gamma}^2+\Lambda\|v\|_{2,\Gamma}^2
\ge \text{min}\{\tilde{l},\Lambda\}\|v\|_{H^1(\Gamma)}^2,$$ so that
the form is coercive whenever $\Lambda>0$. Then, from Lax-Milgram
theorem we get the existence of a solution $v\in H^1(\Gamma)$ of
\eqref{+++}. By the isomorphism property of $-\Delta_\Gamma +1$ is
then follows that $v\in H^2(\Gamma)$, completing the proof of our
claim.

Now, to prove \eqref{Plambda}$_1$, we use \eqref{Plambda}$_2$ in
\eqref{23.2} to get
\begin{equation}\label{28}
\int_\Omega \left (\Delta u-\lambda u+h\right)\Delta
\overline{v}\qquad\text{for all $v\in V$,}
\end{equation}
which clearly implies \eqref{Plambda}$_1$ since for any $\psi\in
L^2(\Omega)$ there are $v\in V$ such that $\Delta v=\psi$, for
example by taking the unique solution $v\in H^2(\Omega)\cap
H^1_0(\Omega)$ with homogeneous Dirichlet boundary conditions. The
proof is then complete.
\end{proof}

The following key estimate shows that the  sesquilinear form
\eqref{16} is appropriate.

\begin{lem}\label{lemma2}
There are positive constants $\lambda_0$ and $C_5$, depending on
$l,k,\Omega,N$, such that for all $\lambda\in \C$, $\text{Re}\lambda
\ge \lambda_0$ we have
$$\text{Re}\,a_\lambda(u,u)\ge C_5 \|u\|_V^2\qquad\text{for all $u\in
V$.}$$
\end{lem}

\begin{proof} By \eqref{16} we have
$$
\begin{aligned}
a_\lambda (u,u)=&\|\Delta u\|_2^2+l \|\Delta_\Gamma
u\|_{2,\Gamma}^2+\lambda \|\nabla u\|_2^2+(\lambda+l)\|d_\Gamma
u\|_{2,\Gamma}^2+\lambda \|u\|_{2,\Gamma}^2\\
-l &\int_\Gamma \Delta u \overline{u}_\nu+k\int _\Gamma u_\nu
\Delta_\Gamma \overline{u}-k\|u_\nu\|_{2,\Gamma}^2-k\int_\Gamma
u_\nu\overline{u}
\end{aligned}
$$
and then
\begin{equation}\label{25}
\begin{aligned}
\text{Re}\,a_\lambda (u,u)=&\|\Delta u\|_2^2+l \|\Delta_\Gamma
u\|_{2,\Gamma}^2+\text{Re}\lambda  \|\nabla
u\|_2^2+(\text{Re}\lambda+l)\|d_\Gamma u\|_{2,\Gamma}^2\\+
&\text{Re}\lambda  \|u\|_{2,\Gamma}^2 +(k-l)\int_\Gamma
\text{Re}[\Delta_\Gamma  u
\overline{u}_\nu]-k\|u_\nu\|_{2,\Gamma}^2-k\int_\Gamma
\text{Re}[u_\nu\overline{u}]\\
\ge &\|\Delta u\|_2^2+l \|\Delta_\Gamma
u\|_{2,\Gamma}^2+\text{Re}\lambda  \|\nabla
u\|_2^2+(\text{Re}\lambda+l)\|d_\Gamma
u\|_{2,\Gamma}^2\\+&\text{Re}\lambda  \|u\|_{2,\Gamma}^2 -(|k|+l)
\int_\Gamma |\Delta_\Gamma
u||u_\nu|-|k|\,\|u_\nu\|_{2,\Gamma}^2-|k|\int_\Gamma \!\!|u_\nu||u|.
\end{aligned}
\end{equation}
By Young inequality we estimate
\begin{equation}\label{26}
|k|\int_\Gamma |u_\nu||u|\le \frac {|k|}2\|u_\nu\|_{2,\Gamma}^2+
\frac {|k|}2\|u\|_{2,\Gamma}^2,
\end{equation}
and, given any $\eps>0$ to be fixed later, by weighted Young
inequality
\begin{equation}\label{27}
(|k|+l)\int_\Gamma |\Delta_\Gamma u||u_\nu |\le \frac {(|k|+l)\eps
}2\|\Delta_\Gamma u\|_{2,\Gamma}^2+ \frac
{|k|+l}{2\eps}\|u_\nu\|_{2,\Gamma}^2.
\end{equation}
Plugging \eqref{26} and \eqref{27} into \eqref{25}, we get
\begin{equation}\label{27.1}
\begin{aligned}
\text{Re}\,a_\lambda (u,u)\ge &\|\Delta u\|_2^2+\left[l-\frac
{(|k|+l)\eps}2 \right] \|\Delta_\Gamma
u\|_{2,\Gamma}^2+\text{Re}\lambda \|\nabla u\|_2^2\\
+&\left(\text{Re}\lambda-\frac {|k|}2 \right)\|u\|_{2,\Gamma}^2
+(\text{Re}\lambda+l)\|d_\Gamma u\|_{2,\Gamma}^2\\-&\left(\frac
{|k|+l}{2\eps}+\frac 32|k|\right)\|u_\nu\|_{2,\Gamma}^2.
\end{aligned}
\end{equation}
Then, by choosing $\eps=l/(|k|+l)$, we get
\begin{equation}\label{29}
\begin{aligned}
\text{Re}\,a_\lambda (u,u)\ge &\|\Delta u\|_2^2+\frac l2
\|\Delta_\Gamma u\|_{2,\Gamma}^2+\text{Re}\lambda \|\nabla u\|_2^2
+\left(\text{Re}\lambda-\frac {|k|}2 \right)\|u\|_{2,\Gamma}^2\\
+&(\text{Re}\lambda+l)\|d_\Gamma
u\|_{2,\Gamma}^2-C_6\|u_\nu\|_{2,\Gamma}^2,
\end{aligned}
\end{equation}
where $C_6=C_6(k,l)=\left(\frac {(|k|+l)^2}{2l}+\frac 32|k|\right)$.
To estimate the last term in the right-hand side of \eqref{29}, we
note that by the embedding $H^{7/4}(\Omega)\hookrightarrow
H^{3/2}(\Omega)$ and by the Trace Theorem there is
$C_7=C_7(\Omega)>0$ such that
$$\|u_\nu\|_{2,\Gamma}^2\le C_7 \|u\|_{H^{7/4}(\Omega)}^2\qquad\text{for all $u\in H^2(\Omega)$.}$$
Consequently, by interpolation inequality (see
\cite{lionsmagenes1}),
$$\|u_\nu\|_{2,\Gamma}^2\le C_7\|u\|_{H^1(\Omega)}^{1/2}\|u\|_{H^2(\Omega)}^{3/2}
\qquad\text{for all $u\in H^2(\Omega)$.}
$$
Using weighted Young inequality we then get, for any $\delta>0$ (to
be fixed below),
$$\|u_\nu\|_{2,\Gamma}^2\le \frac{C_7}{4\delta}\|u\|_{H^1(\Omega)}^2+\frac {3C_7\delta}4\|u\|_{H^2(\Omega)}^2
\qquad\text{for all $u\in H^2(\Omega)$.}
$$
By applying \eqref{10.1} and \eqref{10.2} in the last formula, we get
\begin{equation}\label{30}
\|u_\nu\|_{2,\Gamma}^2\le \frac{C_8}\delta (\|\nabla
u\|_2^2+\|u\|_{2,\Gamma}^2)+C_8\delta \left(\|\Delta
u\|_2^2+\|\Delta_\Gamma u\|_{2,\Gamma}^2+\|u\|_{2,\Gamma}^2\right)
\end{equation}
for all $u\in V$,  where $C_8=C_8(\Omega)>0$. Plugging \eqref{30}
into \eqref{29} we derive
\begin{equation}\label{31}
\begin{aligned}
\text{Re}\,a_\lambda (u,u)\ge &\big(1-C_6C_8\delta\big)\|\Delta
u\|_2^2+\frac l4 \big(2-4C_6C_8\delta/l \big)\|\Delta_\Gamma
u\|_{2,\Gamma}^2\\+&\big(\text{Re}\lambda
-C_6C_8/\delta\big)\|\nabla u\|_2^2 +(\text{Re}\lambda+l)\|d_\Gamma
u\|_{2,\Gamma}^2\\+& \left[\text{Re}\lambda-\tfrac {|k|}2 -C_6C_8
\left(\delta+1/\delta\right)\right]\|u\|_{2,\Gamma}^2.
\end{aligned}
\end{equation}
Choosing $\delta=\delta_0=\text{min }\{2,l\}/(4C_6C_8)$ we rewrite
\eqref{31} as
\begin{equation}\label{32}
\begin{aligned}
\text{Re}\,a_\lambda (u,u)\ge &\frac 12\|\Delta u\|_2^2+\frac l4
\|\Delta_\Gamma u\|_{2,\Gamma}^2+\big(\text{Re}\lambda
-C_6C_8/\delta_0\big)\|\nabla u\|_2^2\\
+&(\text{Re}\lambda+l)\|d_\Gamma u\|_{2,\Gamma}^2+
\left[\text{Re}\lambda-\tfrac {|k|}2 -C_6C_8 \left(\delta_0+
1/\delta_0\right)\right]\|u\|_{2,\Gamma}^2.
\end{aligned}
\end{equation}
Now, by setting
\begin{equation}\label{32.1}
\lambda_0=\text{max }\{C_6C_8/\delta_0, |k|+2C_6C_8
\left(\delta_0+1/\delta_0\right)\},
\end{equation}
we clearly have, when $\text{Re}\lambda\ge \lambda_0$, that
$$\text{Re}\lambda-C_6C_8/\delta_0\ge
0\qquad\text{and }\quad \text{Re}\lambda-\frac{|k|}2-C_6C_8
\left(\delta_0+1/\delta_0\right)\ge\lambda_0/2,$$ so by \eqref{32}
we finally obtain
$$\text{Re}\,a_\lambda (u,u)\ge\frac 12\|\Delta u\|_2^2+\frac l4
\|\Delta_\Gamma u\|_{2,\Gamma}^2+
\frac{\lambda_0}2\|u\|_{2,\Gamma}^2.$$ By setting $C_5=\text{min
}\{\frac 12,\frac l4,\frac {\lambda_0}2\}$ and using \eqref{13} the
proof is complete.
\end{proof}
\begin{rem}\label{remark1}
It is clear from the proof that $\lambda_0\ge 4/l$ and $C_5\le l/4$,
so that $\lambda_0\to+\infty$ and $C_5\to 0$ as $l\to 0^+$. This
instability property will be confirmed in Remark~\ref{remark2}.
\end{rem}
We can now give the
\begin{proof}[\bf Proof of Theorem~\ref{theorem3}.]
By Lemma~\ref{lemma1}, Problem \eqref{Plambda} can be equivalently
written as \eqref{15}. The sesquilinear form $a_\lambda$ in $V$ is
trivially continuous and, by Lemma~\ref{lemma2}, it is also coercive
when $\text{Re}\lambda\ge \lambda_0$. We then apply Lax--Milgram
Theorem (see \cite[p. 376]{dautraylionsvol2}) to get the existence
of a unique solution $u$ of \eqref{Plambda} in $V$. By
Lemma~\ref{lemma1} we also have $u\in H^3$.

We now suppose that $h\in H^s$, $s>1$. To recognize that $u\in
H^{s+2}$ we apply the same bootstrap procedure applied in
Lemma~\ref{lemma1}. More precisely,  we shall prove that for any
$n\in \N_0$ we have
\begin{equation}\label{33}
 u\in H^{\text{min}\{s+2,n+7/2\}}(\Omega),\qquad \text{and}\quad u_{|\Gamma}\in
 H^{\text{min}\{s+2,n+4\}}(\Gamma),
\end{equation}
from which our claim follows for $n$ large enough. We prove
\eqref{33} by induction on $n$. To prove that \eqref{33} holds when
$n=0$ we recognize that, by \eqref{Plambda}$_1$,
$$\Delta u=\lambda u-h\in H^{\text{min}\{s,3\}}(\Omega)\qquad \text{and}\quad u_{|\Gamma}\in
 H^3(\Gamma)$$
so by elliptic regularity (see \cite{lionsmagenes1}) we have $u\in
H^{\text{min}\{s+2,7/2\}}(\Omega)$, which is the required regularity
on $\Omega$ when $n=0$. By the Trace Theorem we then have $u_\nu\in
H^{\text{min}\{s+1/2,2\}}(\Gamma)$. Hence by \eqref{Plambda}$_2$ we
have $-\Delta_\Gamma u+u_{|\Gamma}\in
H^{\text{min}\{s,2\}}(\Gamma)$. By the isomorphism property of
$-\Delta_\Gamma +1$ we then get $u_{|\Gamma}\in
H^{\text{min}\{s+2,4\}}(\Gamma)$ which completes the proof when
$n=0$. To complete the induction process we now suppose that
\eqref{33} holds. Arguing as in the case $n=0$ by
\eqref{Plambda}$_1$ we get
$$\Delta u=\lambda u-h\in H^{\text{min}\{s,n+7/2\}}(\Omega)\qquad \text{and}\quad u_{|\Gamma}\in
 H^{\text{min} \{s+2,4+n\}}(\Gamma)$$
so by elliptic regularity $u\in
H^{\text{min}\{s+2,n+9/2\}}(\Omega)$, By the Trace Theorem we then
have $u_\nu\in H^{\text{min}\{s+1/2,n+3\}}(\Gamma)$, so by using
\eqref{Plambda}$_2$, $-\Delta_\Gamma u+u_{|\Gamma}\in
H^{\text{min}\{s,n+3\}}(\Gamma)$. As before $u_{|\Gamma}\in
H^{\text{min}\{s+2,n+5\}}(\Gamma)$, completing the induction
process.

Finally, to prove \eqref{31xx} we set up the operator $A_\lambda:
D(A_\lambda)\to H^s$, where
\begin{align*}
D(A_\lambda)=&\{(u,v)\in H^{s+2}: (\Delta
u)_{|\Gamma}=ku_\nu+l\Delta_\Gamma v\}\\
\intertext{and} A_\lambda \left(\begin{matrix} u\\
v\end{matrix}\right)=&\left(\begin{matrix} -\Delta u+\lambda u\\
-ku_\nu-l\Delta_\Gamma v+\lambda v\end{matrix}\right).
\end{align*}
One easily sees that $D(A_\lambda)$ is closed in $H^{s+2}$, so it is
an Hilbert space with respect to the scalar product inherited by it.
Moreover $A_\lambda$ is bounded, and $u\in H^{s+2}$ solves
\eqref{Plambda} if and only if $u\in D(A_\lambda)$ and $A_\lambda
u=h$. By previous analysis $A_\lambda$ is bijective, so \eqref{31xx}
follows by the Closed Graph Theorem.
\end{proof}

%%%%%%%%%%%%%%%%%%%%%%%%%%%%%%%%%%%%%%%%%%%%%%%%%%%%%%%%%%%
\section{\bf Analysis of Problem {\bf \eqref{P}}}
\label{section 4}

\noindent We will use here the results of
the  previous section to analyze  Problem \eqref{P}, thus proving
Theorems~\ref{theorem1} and \ref{theorem2}. We start by setting up
the unbounded operator $A:D(A)\subset H\to H$ by
\begin{align}
D(A)=&\{(u,v)\in H^3: (\Delta
u)_{|\Gamma}=ku_\nu+l\Delta_\Gamma v\}\\
\intertext{and}
A \left(\begin{matrix} u\\
v\end{matrix}\right)=&\left(\begin{matrix} \Delta u\\
ku_\nu+l\Delta_\Gamma v\end{matrix}\right).
\end{align}
Our main results are a consequence of the following one.
\begin{thm}\label{theorem4} Operator $A$ generates an analytic
semigroup $\{S(t), t\ge 0\}$ in $H$, and
\begin{equation}\label{38}
\|S(t)\|_{\cal{L}(H)}\le e^{\lambda_0 t}\qquad\text{for all $t\ge
0$},
\end{equation}
where $\lambda_0$ is the positive number given in
Theorem~\ref{theorem3}, so $\{S(t), t\ge 0\}$ is quasi--contractive.
\end{thm}

\begin{proof}
We introduce the unbounded operator $B$ in $H$ by $D(B)=D(A)$ and
$B=A-\lambda_0 I$. Then, given any $u\in D(B)$, we have that $u$
solves \eqref{Plambda} when $\lambda=\lambda_0$ and $h=-Bu$. Hence,
by \eqref{15},
\begin{equation}\label{39}
(Bu,u)_H=-(h,u)_H=-a_{\lambda_0}(u,u)\qquad\text{for all $u\in
D(B)$}.
\end{equation}
Then, by Lemma~\ref{lemma2} we get that $\text{Re}(Bu,u)_H\le 0$,
for all $u\in D(B)$, i.e. $B$ is a dissipative operator in $H$.
Moreover, by Theorem~\ref{theorem3}, $R(I-B)=H$. We then apply
\cite[Theorem~4.6, p.~16]{pazy} to get that $D(B)$ is dense in $H$.
Moreover, given any $u\in D(B)$, by Lemma~\ref{lemma2} and
\eqref{39} we have
\begin{equation}\label{41}
\text{Re}(-Bu,u)_H= \text{Re}a_{\lambda_0}(u,u)\ge C_5 \|u\|_V^2,
\end{equation}
while by \eqref{39} and the continuity of $a_{\lambda_0}$
\begin{equation}\label{42}
|\text{Im }(-Bu,u)_H|\le |a_{\lambda_0}(u,u)|\le C_9 \|u\|_V^2
\end{equation}
for some $C_9=C_9(k,l,N,\Omega)>0$. Combining \eqref{41} and
\eqref{42} we get that $-B$ is a densely defined $m$-sectorial
operator in $H$. We then apply semigroup theory (see for example
\cite[Theorem~5.9, p. 37]{golsteinbook})  which shows that $B$
generates an analytic contraction semigroup $\{T(t), t\ge 0\}$ in
$H$, and consequently $A$ generates an analytic semigroup $\{S(t),
t\ge 0\}$, given by $S(t)=e^{\lambda_0 t}T(t)$, $t\ge 0$, so clearly
\eqref{38} follows.
\end{proof}

Now we can give the proofs of our main results.

\begin{proof}[\bf Proof of Theorem~\ref{theorem1}.]
By Theorem~\ref{theorem4} the operator $A$ generates the analytic,
and hence differentiable, quasi--contractive semigroup $\{S(t), t\ge
0\}$ in $H$. Then, by semigroup theory (see \cite[\S 4.1]{pazy})
given any $u_0\in H$ there is a unique solution
\begin{equation}\label{43}
u\in C([0,\infty);H)\cap C^1((0,\infty);H)
\end{equation}
of the abstract Cauchy problem
\begin{equation}\label{43bis}
\begin{cases}
u'(t)=Au(t),\qquad t>0\\
u(0)=u_0.
\end{cases}
\end{equation}
Clearly, \eqref{43} is nothing but \eqref{33},  and \eqref{43bis} is
the abstract form of Problem \eqref{P}. Moreover, \eqref{1.2} is
nothing but \eqref{38} due to Lemma~\ref{lemma1}. Next, by using the
differentiability property of the semigroup $\{S(t), t\ge 0\}$ and
\cite[\S 2.4]{pazy} we get that $u\in C^\infty((0,\infty);H)$ and
consequently $Bu=Au-\lambda_0u =u'-\lambda_0u \in
C^\infty((0,\infty);H)$. By \eqref{31xx} (when $s=1$) then we get
that $u\in C^\infty((0,\infty);H^3)$. A standard bootstrap procedure
then gives that $u\in C^\infty((0,\infty);H^{2n+1})$ for all
$n\in\N$.  By Morrey's theorem (see for example \cite[Corollaire
IX.13]{brezis} we then get that \eqref{1.3} holds.
\end{proof}

\begin{proof}[\bf Proof of Theorem~\ref{theorem2}.] We introduce, by
recurrence on $n\in\N$, the space
\begin{equation}\label{48}
D(B^n)=\{u\in D(B^{n-1}): Bu\in D(B^{n-1})\}
\end{equation}
endowed with the graph norm
$$\|u\|^2_{D(B^n)}=\sum_{i=0}^n \|D^iu\|^2_H.$$
By Theorem~\ref{theorem3} it is immediate to recognize that
\begin{equation}\label{49}
D(B^n)=\{u\in H^{2n+1}: (\Delta ^i
u)_{|\Gamma}=k(\Delta^{i-1}u)_\nu+l\Delta_\Gamma
(\Delta^{i-1}u)_{|\Gamma}), i=1,\ldots,n\}.
\end{equation}
and that the graph norm it is equivalent to the norm of $H^{2n+1}$
introduced in Section~\ref{section 3}. Since $B$ is a dissipative
operator in $H$ and $R(I-B)=H$ we are able to apply the procedure
outlined in the proof of \cite[Th\'eor\`eme VII.5]{brezis} (see also
\cite [Chapter 1]{breziscazenave}) in the real case, which works as
well in the complex one. Consequently, since $u_0\in D(B^n)$, we get
$$u\in C([0,\infty); D(B^n))\cap
C^1([0,\infty);D(B^{n-1}))\cap\ldots C^n([0,\infty);H)$$ which, by
previous remark, is nothing but \eqref{1.5}.  Finally, if $u_0\in
C^\infty(\Omega)$ and \eqref{1.4} holds for all $i\in\N$ we apply
previous analysis, for any $n\in\N$, together with Morrey's theorem
to get \eqref{1.6}.
\end{proof}

%%%%%%%%%%%%%%%%%%%%%%%%%%%%%%%%%%%%%%%%%%%%%%%%%%%%%%%%%%%
\section{\bf Limit behavior as $\boldsymbol{l\to 0^+}$}
\label{section 5} \noindent This section is devoted to study the
limit behavior of the solution of Problem \eqref{P} when $l\to 0^+$.
The limit problem, at least formally, is given by
\begin{equation}\label{P0}
\begin{cases} u_t-\Delta u=0 \qquad &\text{in
$Q=(0,\infty)\times\Omega$,}\\
 u_t=ku_\nu \qquad &\text{on
$[0,\infty)\times\Gamma$,}\\
 u(0,x)=u_0(x) &
 \text{in $ \Omega$,}
\end{cases}\end{equation}
which has been studied in \cite{vazvit} (see also \cite{bandle}). We recall the following definition and result from that paper. In what follows we restrict to the real-valued case.

\begin{definition}[{\cite[\bf Definition~1]{vazvit}}] Given $u_0\in H^1(\Omega)$ we say that
\begin{equation}\label{5.1}
u\in C([0,T);H^1(\Omega)),\qquad T>0,
\end{equation}
is a weak solution of \eqref{P0} if
\begin{equation}\label{5.2}
u(0)=u_0
\end{equation}
and
\begin{equation}\label{5.3}
-\int_0^T\int_\Omega u\varphi_t+\int_0^T\int_\Omega \nabla u\nabla
\varphi+\frac 1k \int_0^T\int_\Gamma u\varphi_t=0
\end{equation}
for all $\varphi\in C^\infty_c((0,T)\times \overline{\Omega})$.
\end{definition}

\begin{thm}[{\cite[\bf First part of Theorem~1]{vazvit}}]\label{theorem5} If $N\ge 2$  there is $u_0\in
C^\infty(\overline{\Omega})$ satisfying the compatibility conditions
$$
\Delta ^{n}u_0=k(\Delta^{n-1} u_0)_\nu \qquad\text{on $\Gamma$ for
all $n\in\N$,}
$$
such that Problem \eqref{P0} has no weak solutions.
\end{thm}

The first step in our analysis is the following
\begin{lem}\label{lemma3} Theorem~\ref{theorem5} holds also if
\eqref{5.1} is weakened to
\begin{equation}\label{5.4}
 u\in C_w([0,T);H^1(\Omega)),
\end{equation}
that is it concerns also weakly continuous solutions.
\end{lem}
\begin{proof}
Looking at the proof in the quoted paper one immediately sees that
the continuity of $u$ was used only at two places: at first in order
that \eqref{5.3} makes sense, and at second to recognize that the
functions $t\mapsto \langle u(t),\Phi'_n\rangle$ are continuous in
$[0,T)$, where $\langle\cdot,\cdot\rangle$ denotes an equivalent
scalar product in $H^1(\Omega)$ and $\Phi'_n$, $n\in\N$ are
eigenfunctions of a suitable eigenvalue problem, which belong to
$C^\infty(\overline{\Omega})$. Both facts continue to hold when
\eqref{5.1} is weakened to \eqref{5.4}.
\end{proof}

We can now state the main result of this section.
\begin{thm}\label{theorem6} Let $u_0\in H^1(\Omega)$ be an initial datum such that
Problem \eqref{P0} has no weak solutions $u\in
C_w([0,T);H^1(\Omega))$ for any $T>0$, and denote by $u^l$ the
solution of \eqref{P} corresponding to $u_0$ and $l$ given by
Theorem~\ref{theorem1}. \begin{footnote}{which is real valued since
$u_0$ is real valued}\end{footnote} Then, for any $T>0$, we have
\begin{equation}\label{5.5}
\max_{t\in [0,T]} \|u^l(t)\|_{H^1(\Omega)}\to\infty\qquad\text{as \
$l\to 0^+$.}
\end{equation}
\end{thm}

\begin{proof} We suppose by contradiction that \eqref{5.5} fails for
some fixed $T>0$. Then there is a sequence $l_n\to 0^+$ such that,
\begin{equation}\label{5.6}
\|u^n\|_{C([0,T];H^1(\Omega))}\le C_{10}\qquad \text{for all
$n\in\N$,}
\end{equation}
were we denoted $u^n=u^{l_n}$ for simplicity, and
$C_{10}=C_{10}(T,u_0,\Omega)>0$. Since, by Theorem~\ref{theorem1},
we have $u^n\in C^\infty((0,\infty)\times\overline{\Omega})$ we are
allowed, for any $t\in(0,T)$, to multiply the heat equation by a
test function $v\in C^\infty_c(\Omega)$ and  to integrate by parts
to get
\begin{equation}\label{5.6bis}
\left|\int_\Omega u^n_t(t)v\right|\le
\|u^n(t)\|_{H^1(\Omega)}\|\nabla v\|_2.
\end{equation}
By a standard density argument \eqref{5.6bis} holds true for all
$v\in H^1_0(\Omega)$ and then, by \eqref{5.6}, we get the second
estimate we need, that is
\begin{equation}\label{5.7}
\|u^n_t\|_{L^\infty((0,T);H^{-1}(\Omega))}\le C_{10}\qquad \text{for
all $n\in\N$.}
\end{equation}
By \eqref{5.6} we get that, up to a subsequence,
\begin{equation}\label{5.8}
u^n\to u\qquad\text{weakly$^*$ in $L^\infty((0,T);H^1(\Omega))$}.
\end{equation}
Moreover, by combining \eqref{5.6} with \eqref{5.7} and using the
compactness of the embedding $H^1(\Omega)\hookrightarrow
L^2(\Omega)\hookrightarrow H^{-1}(\Omega)$ we also  get, by
Aubin--Lions compactness lemma, that
\begin{equation}\label{5.9}
u^n\to v\qquad\text{strongly in $C([0,T];H^{-1}(\Omega))$}.
\end{equation}
Moreover, by \eqref{5.8} we get that $u^n\to u$ weakly$^*$ in
$L^\infty((0,T);H^{-1}(\Omega))$, while by \eqref{5.9} $u^n\to v$ in
the same sense, so $u=v$. Then we can combine
\eqref{5.8}--\eqref{5.9} to
\begin{equation}\label{5.10}
u^n\to u\quad\text{weakly$^*$ in $L^\infty((0,T);H^1(\Omega))$ and
strongly in $C([0,T];H^{-1}(\Omega))$}.
\end{equation}

We now claim that $u$ is a weak solution of \eqref{P0} in the class
\eqref{5.4}. Once this claim is proved the proof is complete, since
we are in contradiction with Theorem~\ref{theorem5} as extended by
Lemma~\ref{lemma3}.  By applying \cite[Theorem~2.1]{strauss}, we get
\begin{equation}\label{5.11}
u\in C_w([0,T];H^1(\Omega)).
\end{equation}
Moreover, by \eqref{5.10} it immediately follows that \eqref{5.2}
holds. Multiplying the heat equation in \eqref{P} by a test function
$\psi\in C^\infty_c(0,T)\times\overline{\Omega}$, integrating by
parts in $\Omega$, using the boundary condition in \eqref{P}, using
\eqref{1} and finally integrating by parts in time twice we get that
$u^n$ satisfies the distribution identity
\begin{equation}\label{5.12}
-\int_0^T\int_\Omega u^n\varphi_t+\int_0^T\int_\Omega \nabla
u^n\nabla \varphi+\frac 1k \int_0^T\int_\Gamma u^n\varphi_t-\frac
{l_n}k\int_0^T\int_\Gamma u^n\Delta_\Gamma \varphi=0
\end{equation}
for all $\varphi\in C^\infty_c((0,T)\times \overline{\Omega})$. By
\eqref{5.10} we can pass to the limit as $n\to\infty$ in
\eqref{5.12}, and we get that $u$ satisfies \eqref{5.3}, completing
the proof of our claim.
\end{proof}

\begin{rem}\label{remark2}
Theorem~\ref{theorem6} shows that the instability property of
$\lambda_0=\lambda_0(l)$ pointed out in Remark~\ref{remark1} does
not depend on our estimates. Indeed, if one could improve them to get
a sequence $l_n\to 0^+$ and
$\lambda_0(l_n)\le\overline{\lambda}<\infty$, then by \eqref{1.2}
one would contradict \eqref{5.5}.
\end{rem}

%%%%%%%%%%%%%%%%%%%%%%%%%%%%%%%%%%%%%%%%%%%%%%%%%%%%%%%%%%%
\section{\bf Open problems and final remarks}
\label{section 6} \noindent Although
Theorems~\ref{theorem1}--\ref{theorem2} give existence and
uniqueness of solutions to Problem \eqref{P} in a Hilbert
framework, building a satisfactory theory for $C^\infty$ initial
data, many interesting problems are still open, both of theoretical
and applied nature.
\renewcommand{\labelenumi}{\bf \arabic{enumi}.}
\begin{enumerate}
\item We are not able to produce a satisfactory regularity theory in
even order spaces $H^{2n}$, $n\ge 1$, which is particularly bad for
$n=1$. A new estimate in $V$ would be necessary.
\item The extension of the analysis to more general problems, like
the ones considered by \cite{coclite1}--\cite{coclite2} has still to
be done. In particular, Lemma~\ref{lemma1} has to be properly
extended.
\item Our arguments, which are based on Lax--Milgram theorem, cannot
be extended to the case of Banach spaces. Now, it would be natural
to consider the case of $u_0\in W^{1,p}(\Omega)$,
${u_0}_{|\Gamma}\in W^{1,p}(\Gamma)$.
\item We do not give explicit representation formulas of the
solution $u$, even for regular data. In particular we are not able
to apply the Fourier method, which would be based on the study of the
eigenvalue problem
\begin{equation}\label{Plambdabis}
\begin{cases} -\Delta u=\lambda u \qquad &\text{in
$\Omega$,}\\
-ku_\nu-l\Delta_\Gamma u=\lambda u\qquad &\text{on $\Gamma$.}
\end{cases}\end{equation}
The elliptic theory developed in Section~\ref{section 3} allows to
prove in a simple way that
$$\Sigma:=\{\lambda\in \C: \text{\eqref{Plambdabis} has a nontrivial solution
$u\in V$}\}$$ is at most countable, but it is far from   giving an
exhaustive spectral theory, since the operator $A$ or equivalently
$A_\lambda$ is not symmetric in $H$. Actually formula \eqref{15bis}
suggest suggest some symmetry of  the operator $A_\lambda$, but in a
framework of Krein spaces. This study is left to specialists in
Krein spaces theory.
\end{enumerate}

\vskip .5cm

%%%%%%%%%%%%%%%%%%%%%%%%%%%%%%%%%%%%%%%%%%%%%%%%%%%%%%%%%%
\noindent \textsc{Acknowledgment.}  The first author was partially
supported by Spanish Project MTM2008-06326-C02  and by ESF Programme
``Global and geometric aspects of nonlinear partial differential
equations". The second author was partially supported by the
M.I.U.R. Project "Metodi variazionali ed equazioni differenziali
nonlineari" (Italy). This author wishes to thank Univ.~Aut\'onoma de
Madrid for its kind hospitality.

\vskip .5cm
\newpage

%\bibliographystyle{amsplain}
%\bibliography{biblio2}

\begin{thebibliography}{10}

\bibitem{AMP}
W.~Arendt, G.~Metafune, D.~Pallara, and S.~Romanelli, \emph{The
{L}aplacian
  with {W}entzell-{R}obin boundary conditions on spaces of continuous
  functions}, Semigroup Forum \textbf{67} (2003), no.~2, 247--261.

\bibitem{bandlereichel}
C.~Bandle and W.~Reichel, \emph{A linear parabolic problem with
  non--dissipative dynamical boundary conditions}, Proceedings of the 2004
  Swiss--Japanese Seminar (Zurich, 6--10 December 2004) (Michel Chipot and
  Hirokazu Ninomiya, eds.), World Scientific, March 2006, pp.~45--77.

\bibitem{bandle}
C.~Bandle, J.~von Below, and W.~Reichel, \emph{Parabolic problems
with
  dynamical boundary conditions: eigenvalue expansions and blow up}, Atti
  Accad. Naz. Lincei Cl. Sci. Fis. Mat. Natur. Rend. Lincei (9) Mat. Appl.
  \textbf{17} (2006), no.~1, 35--67.

\bibitem{BHLN}
P.~Binding, R.~Hryniv, H.~Langer, and B.~Najman, \emph{Elliptic
eigenvalue
  problems with eigenparameter dependent boundary conditions}, J. Differential
  Equations \textbf{174} (2001), no.~1, 30--54.

\bibitem{brezis}
H.~Brezis, \emph{Analyse fonctionnelle}, Masson, Paris, 1983,
Th\'eorie et
  applications.

\bibitem{breziscazenave}
H.~Brezis and T.~Cazenave, \emph{Nonlinear evolution equations,
Unpublished}.

\bibitem{carmugnolo1}
S.~Cardanobile and D.~Mugnolo, \emph{Qualitative properties of
coupled
  parabolic systems of evolution equations}, Ann. Sc. Norm. Super. Pisa Cl.
  Sci. (5) \textbf{7} (2008), no.~2, 287--312.

\bibitem{carmugnolo2}
\bysame, \emph{Parabolic systems with coupled boundary conditions},
J.
  Differential Equations \textbf{247} (2009), no.~4, 1229--1248.

\bibitem{coclite1}
G.~M. Coclite, A.~Favini, G.~Ruiz Goldstein, J.~A. Goldstein, and
S.~Romanelli,
  \emph{Continuous dependence on the boundary conditions for the {W}entzell
  {L}aplacian}, Semigroup Forum \textbf{77} (2008), no.~1, 101--108.

\bibitem{coclite2}
G.~M. Coclite, G.~R. Goldstein, and J.~A. Goldstein, \emph{Stability
estimates
  for parabolic problems with {W}entzell boundary conditions}, J. Differential
  Equations \textbf{245} (2008), no.~9, 2595--2626.

\bibitem{dautraylionsvol2}
R.~Dautray and J.-L. Lions, \emph{Mathematical analysis and
numerical methods
  for science and technology. {V}ol. 2}, Springer-Verlag, Berlin, 1988,
  Functional and variational methods.

\bibitem{Eng}
K.~J. Engel, \emph{The {L}aplacian on {$C(\overline\Omega)$} with
generalized
  {W}entzell boundary conditions}, Arch. Math. (Basel) \textbf{81} (2003),
  no.~5, 548--558.

\bibitem{ES1}
J.~Ercolano and M.~Schechter, \emph{Spectral theory for operators
generated by
  elliptic boundary problems with eigenvalue parameter in boundary conditions.
  {I}}, Comm. Pure Appl. Math. \textbf{18} (1965), 83--105.

\bibitem{ES2}
\bysame, \emph{Spectral theory for operators generated by elliptic
boundary
  problems with eigenvalue parameter in boundary conditions. {II}}, Comm. Pure
  Appl. Math. \textbf{18} (1965), 397--414.

\bibitem{escher1}
J.~Escher, \emph{Global existence and nonexistence for semilinear
parabolic
  systems with nonlinear boundary conditions}, Math. Ann. \textbf{284} (1989),
  285--305.

\bibitem{escher3}
\bysame, \emph{A note on quasilinear parabolic systems with
dynamical boundary
  conditions}, Progress in partial differential equations: the Metz surveys, 2
  (1992), Longman Sci. Tech., Harlow, 1993, pp.~138--148.

\bibitem{escher2}
\bysame, \emph{Quasilinear parabolic systems with dynamical boundary
  conditions}, Comm. Partial Differential Equations \textbf{18} (1993),
  1309--1364.

\bibitem{escher4}
\bysame, \emph{On the qualitative behaviour of some semilinear
parabolic
  problems}, Differential Integral Equations \textbf{8} (1995), no.~2,
  247--267.

\bibitem{fgg}
A.~Favini, G.~R. Goldstein, J.~A. Goldstein, and S.~Romanelli,
\emph{The heat
  equation with generalized {W}entzell boundary condition}, J. Evol. Equ.
  \textbf{2} (2002), no.~1, 1--19.

\bibitem{FGGRLLB}
A.~Favini, G.~R.Goldstein, J.A. Goldstein, and S~Romanelli,
\emph{The heat
  equation with dynamic boundary conditions}, in preparation.

\bibitem{golsteinbook}
J.~A. Goldstein, \emph{Semigroups of linear operators and
applications}, Oxford
  Mathematical Monographs, The Clarendon Press Oxford University Press, New
  York, 1985.

\bibitem{goldsteingiselle}
G.~Goldstein~R., \emph{Derivation and physical interpretation of
general
  boundary conditions}, Adv. Differential Equations \textbf{11} (2006), no.~4,
  457--480.

\bibitem{grobbelaar}
M.~Grobbelaar Van~Dalsen, \emph{Semilinear evolution equations and
fractional
  powers of a closed pair of operators}, Proc. Roy. Soc. Edinburgh Sect. A
  \textbf{105} (1987), 101--115.

\bibitem{hintermann}
T.~Hintermann, \emph{Evolution equations with dynamic boundary
conditions},
  Proc. Roy. Soc. Edinburgh Sect. A \textbf{113} (1989), no.~1-2, 43--60.

\bibitem{jost}
J.~Jost, \emph{Riemannian geometry and geometric analysis}, fifth
ed.,
  Universitext, Springer-Verlag, Berlin, 2008.

\bibitem{lionsmagenes1}
J.-L. Lions and E.~Magenes, \emph{Probl\`emes aux limites non
homog\`enes et
  applications. {V}ol. 1}, Travaux et Recherches Math\'ematiques, No. 17,
  Dunod, Paris, 1968.

\bibitem{mugnoloarxiv1}
D.~Mugnolo, \emph{Vector--valued heat equations and networks with
coupled
  dynamic boundary conditions},
  \verb+http://arxiv.org/PS_cache/arxiv/pdf/0903/0903.3580v3.pdf+.

\bibitem{pazy}
A.~Pazy, \emph{Semigroups of linear operators and applications to
partial
  differential equations}, Applied Mathematical Sciences, vol.~44,
  Springer-Verlag, New York, 1983.

\bibitem{strauss}
W.~A. Strauss, \emph{On continuity of functions with values in
various {B}anach
  spaces}, Pacific J. Math. \textbf{19} (1966), 543--551.

\bibitem{taylor}
M.~E. Taylor, \emph{Partial differential equations}, Texts in
Applied
  Mathematics, vol.~23, Springer-Verlag, New York, 1996, Basic theory.

\bibitem{vazvit2}
J.~L. Vazquez and E.~Vitillaro, \emph{Heat equation with dynamical
boundary
  conditions of locally reactive type}, Semigroup Forum \textbf{74} (2007),
  no.~1, 1--40.

\bibitem{vazvit}
\bysame, \emph{Heat equation with dynamical boundary conditions of
reactive
  type}, Comm. Partial Differential Equations \textbf{33} (2008), no.~4-6,
  561--612.

\bibitem{vazvitLLB}
\bysame, \emph{On the {L}aplace equation with dynamical boundary
conditions of
  reactive-diffusive type}, J. Math. Anal. Appl. \textbf{354} (2009), no.~2,
  674--688.

\bibitem{ziemer}
W.~P. Ziemer, \emph{Weakly differentiable functions}, Graduate Texts
in
  Mathematics, vol. 120, Springer-Verlag, New York, 1989.

\end{thebibliography}

\def\cprime{$'$}
\providecommand{\bysame}{\leavevmode\hbox
to3em{\hrulefill}\thinspace}
\providecommand{\MR}{\relax\ifhmode\unskip\space\fi MR }
% \MRhref is called by the amsart/book/proc definition of \MR.
\providecommand{\MRhref}[2]{%
  \href{http://www.ams.org/mathscinet-getitem?mr=#1}{#2}
} \providecommand{\href}[2]{#2}

\end{document}